\documentclass[a4paper, 
10pt
]{article}

\usepackage{amsmath}
\usepackage{multirow}
\usepackage[table]{xcolor}

\usepackage{booktabs}

\usepackage[table]{xcolor}

\usepackage{hyperref}
\usepackage{siunitx}

\newcommand{\mupar}{\ensuremath{\boldsymbol{\mu}}}

\newcommand{\vtilde}{\ensuremath{\tilde{\boldsymbol{v}}^{\text{wake}}}}
\newcommand{\vhat}{\ensuremath{\hat{\boldsymbol{v}}^{\text{wake}}}}
\newcommand{\vtarget}{\ensuremath{\boldsymbol{v}^{\text{target}}}}
\newcommand{\vwake}{\ensuremath{\boldsymbol{v}^{\text{wake}}}}

\newcommand{\nn}{\ensuremath{\text{ANN}}}
\newcommand{\autoe}{\ensuremath{\text{AE}}}
\newcommand{\rom}{\ensuremath{\text{ROM}}}
\newcommand{\podec}{\ensuremath{\text{POD}}}
\newcommand{\rbf}{\ensuremath{\text{RBF}}}
\newcommand{\ga}{\ensuremath{\text{GA}}}
\newcommand{\cfd}{\ensuremath{\text{CFD}}}

\usepackage{cases}
\definecolor{mediumpurple}{rgb}{0.58, 0.0, 0.83}
\definecolor{black}{rgb}{0.0, 0.0, 0.0}
\newcommand{\RA}[1]{{\color{black}#1}}
\newcommand{\RB}[1]{{\color{black}#1}}
\newcommand{\RC}[1]{{\color{black}#1}}
\newcommand{\RALL}[1]{{\color{black}#1}}

\usepackage{subfig}
\usepackage{array,multirow}

\usepackage{authblk}
\usepackage[a4paper,left=3cm,right=3cm,top=2.5cm,bottom=2.5cm]{geometry}
\usepackage{graphicx}
\usepackage{amssymb}
\newtheorem{remark}{Remark}%
\newcommand{\bmhead}{\paragraph}
\date{}

\begin{document}

\title{Towards a machine learning pipeline in reduced order modelling for inverse problems: neural networks for boundary parametrization, dimensionality reduction and solution manifold approximation}

\author[]{Anna Ivagnes\footnote{anna.ivagnes@sissa.it}}
\author[]{Nicola Demo\footnote{nicola.demo@sissa.it}}
\author[]{Gianluigi Rozza\footnote{gianluigi.rozza@sissa.it}}

\affil{Mathematics Area, mathLab, SISSA, via Bonomea 265, I-34136
  Trieste, Italy}

\maketitle

\begin{abstract}
In this work, we propose a model order reduction framework to deal with 
inverse problems in a non-intrusive setting. Inverse problems, especially in a partial differential equation context, require a huge computational load due to the iterative optimization process. To accelerate such a procedure, we apply a numerical pipeline that involves artificial neural networks to parametrize the boundary conditions of the problem in hand, compress the dimensionality of the (full-order) snapshots, and approximate the parametric solution manifold. It derives a general framework capable to provide an ad-hoc parametrization of the inlet boundary and quickly converges to the optimal solution thanks to model order reduction. We present in this contribution the results obtained by applying such methods to two different CFD test cases.
\end{abstract}


\section{Introduction}
\label{sec:intro}
Inverse problems is a wide family of problems that crosses many different sciences and engineering fields. Inverse problems aim to compute from the given observations the cause that has produced them, as also explained in \cite{cetrangolo2021reduced,richter2020inverse}. Formally, we can
consider an input $I$ and an output $O$, and suppose that there exists a map
\[\mathcal{M}:i \rightarrow o \]
that models a mathematical or physical law. The computation of the output as $o=\mathcal{M}(i)$ is called \emph{direct problem}, whereas the problem of finding the input given the output is called \emph{inverse problem}.
Given a certain output $o_t$, the inverse problem consists of inverting the map $\mathcal{M}$ and finding the input $i_t$ which produces the output $o_t$, i.e., which satisfies $\mathcal{M}(i_t)=o_t$. 
\RB{
Inverse problems may be of interest for a lot of mathematical fields, from heat transfer problems to fluid dynamics. The case study which is here analysed is a Navier-Stokes flow in a circular cylinder, and the aim is to find the proper boundary condition in order to obtain the expected distribution within the domain. Such an application tries to solve a typical naval problem, numerically looking for the inlet setting which provides the right working condition during the optimization of the propulsion system. Propeller optimization is indeed very sensitive to modifications in the velocity distribution at the propeller plane: to obtain an optimized artifact it becomes very important to set up the numeric simulation such that the velocity distribution is as close as possible to the target distribution, usually collected by experimental tests.
The problem is the identification of the inlet boundary condition, given the velocity distribution at the so-called \emph{wake}, which is the plane (or a portion of it) orthogonal to the cylinder axis where the propeller operates. To achieve that, the inlet distribution is searched by parametrizing the target wake --- by exploiting a neural network, as we will describe in the next paragraphs --- and optimizing these parameters such that in the simulation the velocity is close to the target wake.
It must be said that to produce meaningful results, we assume here the flow has a main direction that is perpendicular to the inlet and wake planes: in such a way, the distributions at these planes are similar to each other, allowing us to search for the optimal inlet among the parametrized wake distributions.
Even if in this case the numerical experiments are conducted in a Computational Fluid Dynamics (CFD) setting, the methodology is in principle easily transferable to different contexts.
}

Typically, such an operation is performed within an optimization procedure, in which the direct problem is iteratively solved by varying the input until the desired output is reached. This, of course, implies the necessity to numerically parametrize the input in a proper way, possibly allowing a large variety of admissible input\RC{s} and at the same time a limited number of parameters.
Moreover, the necessity to solve the direct problem for many different instances makes the entire process computationally expensive, especially dealing with the numerical solution of Partial Differential Equations (PDEs). A possible solution to overcome such computational burden is offered by the Reduced Order Modelling (\rom) techniques.\\

\rom{} constitutes a constantly growing approach for model simplification, allowing for a real-time approximation of the numerical solution of the problem at hand. Among the methods already introduced in the literature, the Proper Orthogonal Decomposition (\podec) has become in recent developments an important tool for dealing with PDEs, especially in parametric settings~\cite{salmoiraghi2018free,quarteroni2015reduced,rozza2015book,morhandbook2020}. Such a framework aims to efficiently combine the numerical solutions for different configurations of the problem, typically pre-computed using consolidated methods --- e.g. finite volume, finite element --- such that at any model inference all this information is combined for providing a fast approximation. Within iterative and many-query processes, like inverse problems, this methodology allows a huge computational gain. The many repetitions of the direct problem, needed to find the target input, can be performed at the reduced level, requiring then a finite and fixed set of numerical solutions only for building the \rom{}. The coupling between \rom{} and the inverse problem has been already investigated in literature for heat flux estimation in a data assimilation context~\cite{Morelli1}, in aerodynamic application~\cite{bui2004aerodynamic}, in h\RC{a}emodynamic problems~\cite{lassila2013}. An alternative way to efficiently deal with this kind of problem has been explored in a Bayesian framework~\cite{li2014adaptive}. \RC{Moreover,} among all the contributions in literature we cite~\cite{VITALE2012788,HuangInverseProblem} as an example of inverse problem with pointwise observations and inverse problem in a boundary element method context, respectively.\\

\begin{figure}[ht]
    \centering
    \includegraphics[width=\textwidth]{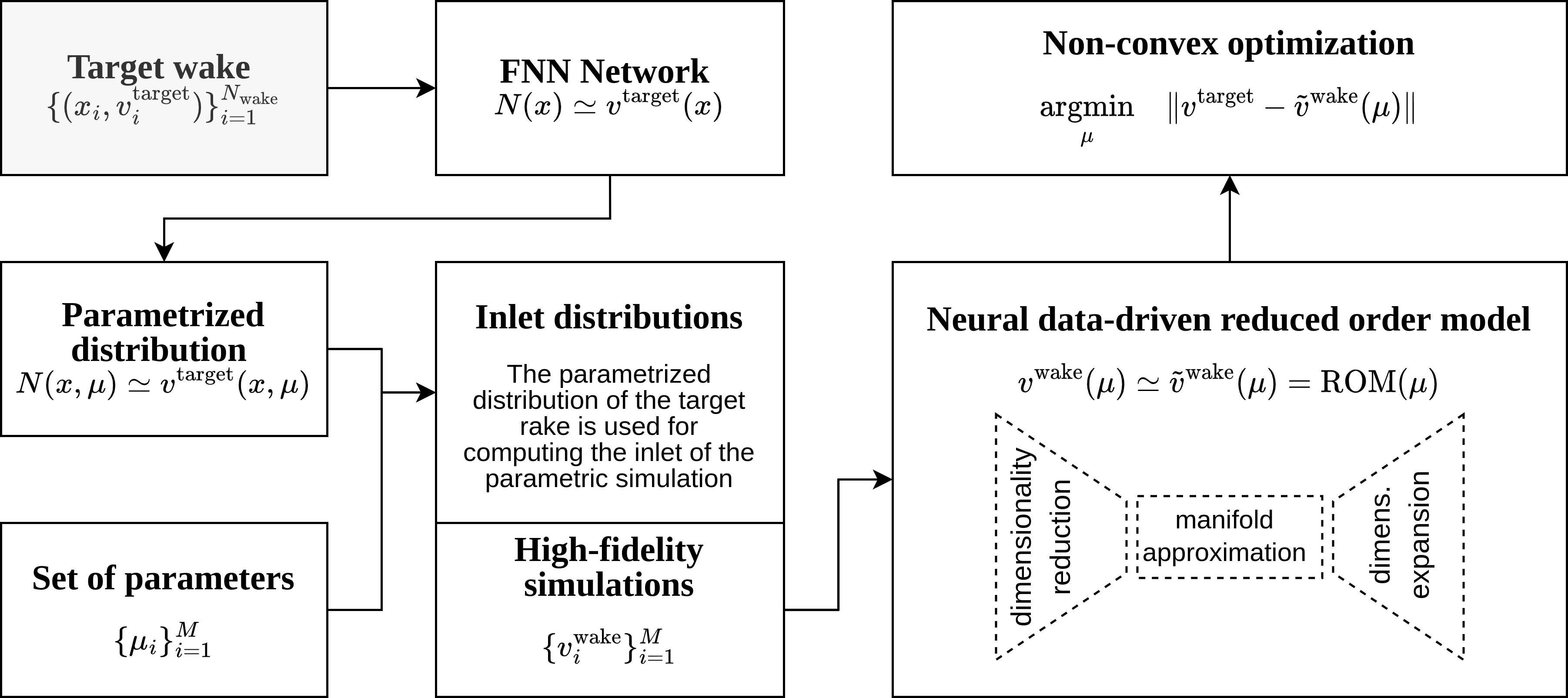}
    \caption{Flow diagram for the data-driven pipeline followed in the paper.}
    \label{fig:diagram}
\end{figure}

This contribution introduces an entire and novel machine learning pipeline to deal with the inverse problems in a \rom{} setting. In specific, we combine three different use\RC{s} of Artificial Neural Network (\nn{}), that are: {\it i}) parametrization of the boundary condition given a certain target distribution or pointwise observations, {\it ii}) dimensionality compression of the discrete space of the original --- the so-called full-order --- model and {\it iii}) approximation of the parametric solution manifold. It derives a data-driven pipeline (graphically represented in Figure~\ref{fig:diagram}) able to provide a parametrization of the original problem, which is successively exploited for the optimization in the reduced space. Finally, the optimization is carried out by involving a Genetic Algorithm (GA), but in principle can be substituted by any other optimization algorithm.

The contribution presents in Section~\ref{sec:other} an algorithmic overview of the employed methods, whereas Section~\ref{sec:results} illustrates the results of the numerical investigation pursued to the above-mentioned test case. In particular, we present details for all the intermediate outcomes, comparing them to the results obtained by employing state-of-the-art techniques. Finally, Section~\ref{sec:conclusions} is dedicated to summarizing the entire content of the contribution,
drawing future perspectives and highlighting the criticisms highlighted during the development of this contribution.

\section{Methodology}
\label{sec:other}
We dedicate this section to providing an algorithmic overview of the numerical tools composing the computational pipeline.
\newcommand{\ANN}{\text{ANN}}
\label{sec:methods}

\subsection{Boundary parametrization using \ANN}
\label{subsec:ANN_param}
\RC{Neural networks are a class of regression techniques and the general category of \emph{Feed-forward} neural networks has been the subject of considerable research in several fields in recent years.} The capability of \ANN{} to approximate any function~\cite{hornik_1989} and the even greater computational availability allowed indeed the massive employment of such an approach to overcome many limitations.
In the field of PDEs, we cite~\cite{raissi2019physics,prnn_paper,pod_nn_paper,lee2020model,pichi2021artificial,Lu2021,wang2021learning,PapapiccoDemoGirfoglioStabileRozza2021} as some of main impacting frameworks recently explored. A not yet investigated usage of \ANN{}, to the best of authors' knowledge, is instead the parametrization of a given (scattered) function.
\RC{We suppose that we have a target distribution
$ \vtarget = (\vtarget_i)_{i=1}^{P}$, corresponding to the wake velocity distribution in our case, which has $P$ degrees of freedom. We want to reproduce this distribution by exploiting a neural network technique.}\\

\RC{A neural network is defined as a concatenation of an input layer, multiple hidden layers, and a final layer of output neurons. An output of the $i$-th neuron in the $l$-th layer of the network is generally defined as:

\begin{equation}
    h^l_i = \sigma \left( \sum_{j=1}^{N_{l-1}} W^l_{ij} h^{l-1}_j +b^l_i \right), \, i=1, \dots, N_l \, ,l=1, \dots,H \,.
    \label{eq:nn}
\end{equation}
}
where $\sigma$ is the activation function (which provides non-linearity), $b\RC{^l_i}$ the bias and $W$ refers to the so-called weights of the synapses of the network, \RC{$N_l$ is the number of neurons of the $l$-th hidden layer, $H$ is the number of layers of the network.}

The bias and the weights are the hyperparameters of the network, that are tuned during the training procedure to converge to the optimal values for the approximation in hand. 
We can think \RC{of} these hyperparameters as the degree of freedom of our approximation, allowing us to manipulate such distribution by perturbating them. We define the optimal hyperparameters (computed with a generic optimization procedure \cite{rumelhart1986learning, rojas1996backpropagation}) as $b^*$ and $W^*$; \RC{the network exploiting these hyperparameters reproduces as output an approximation of our target wake distribution:
\[ N(\mathbf{x}) \simeq \vtarget(\mathbf{x}) .\]}
The input $\mathbf{x}$ \RC{to the whole neural network} in this paper corresponds to the polar coordinates of the points of the wake, so we have a two-dimensional input.
\RC{We can then} rearrange Eq.~\ref{eq:nn} to \RC{express the parametrized output of a single hidden layer as follows:
\begin{equation}
    h^l_i(\mupar) = \sigma \left( \sum_{j=1}^{N_{l-1}} W^{*\,l}_{ij} h^{l-1}_j + (b^{* \, l}_i + \mu_i^l) \right), \, i=1, \dots, N_l \, ,l=1, \dots,H.
    \label{eq:nn2}
\end{equation}
where \RC{$\mupar^l$} is the vector of parameters in layer $l$, which applies only to the bias of the layers.
We finally obtain the parametrized output $N(\mathbf{x}, \mupar)$.}

The main advantage of this approach is the capability to parametrize any initial distribution, since the weights are initially learnt during the training and then manipulated to obtain its parametric version. On the other hand, the dimension of the weights \RC{is} typically very large, especially in networks with more than one layer. 
\RC{In networks composed of a large number of layers, a high number of hyperparameters should be tuned.}
A possible solution to overcome such a\RC{n} issue could be to manipulate just a subset of all the hyperparameters. Indeed, in this paper, only the bias parameters of two hidden layers are perturbed.

Such a posteriori parametrization of any generic distribution is employed in this work to deal with the inverse problem. The main idea is to compute different inlet velocity distributions corresponding to different weights of the \ANN{}. 
The weights perturbations are used as parameters to build the reduced order model, providing an \emph{ad-hoc} parametrization of the boundary condition based on the target output. It must be said that such parametrization may lead \RC{to} good results only if some correlation between the boundaries and the target output exists.

\subsection{Model Order Reduction}
\rom{}s are a class of techniques aimed to reduce the computational complexity of mathematical models.

The technique used in this paper is \emph{data-driven} or \emph{non-intrusive} ROM, which allows \RC{us to} build a reduced model without requiring the knowledge of the governing equations of the observed phenomenon. For this reason, this technique is suitable to deal with experimental data and it has been widely used in numerical simulations for industrial, naval, biomedical, and environmental applications~\cite{tezzele2018ecmi,demoortaligustinrozzalavini2020bumi,tezzele2020enhancing,georgaka2020hybrid,dutta2021greedy,aria2021,Girfoglio2020_b}. 

Non-intrusive \rom{}s are based on an offline-online procedure. Each stage of the procedure can be approached in different ways, which are analyzed in the following Sections. All the techniques presented in the next lines have been implemented in the Python package called  EZyRB~\cite{DemoTezzeleRozza2018EZyRB}.

\subsubsection{Dimensionality reduction}
\RB{
In the dimensionality reduction step, a given set of high-dimensional snapshots is represented by a few reduced coordinates, in order to fight the curse of dimensionality and make the pipeline more efficient.
}

We consider the following matrix of $M$ snapshots, collecting our data:
\[
\begin{split}
\mathbf{Y}=
&\begin{bmatrix}
\vert & \vert & \vert & &\vert\\
\vwake_1&\vwake_2&\vwake_3&...&\vwake_M\\
\vert&\vert&\vert& &\vert
\end{bmatrix}= \\
&\begin{bmatrix}
\vert & \vert & \vert & &\vert\\
\vwake(\mu_1)&\vwake(\mu_2)&\vwake(\mu_3)&\dotsc&\vwake(\mu_M)\\
\vert&\vert&\vert& &\vert
\end{bmatrix},
\end{split}
\]
where $\vwake_i \in \mathbb{R}^P$, $i=1,\dotsc,M$ are the velocity distributions in our case, each one corresponding to a different set of parameters $\mu_i \in \mathbb{R}^p$ for $i=1,\dotsc,M$. 
\RB{The goal is to calculate the \emph{reduced} snapshots $\{ \vhat_i \in \mathbb R^L\}_{i=1}^M$ such that
\begin{equation}
\vtilde_i~\simeq~\Psi^{-1}(\Psi(\vwake_i)), \quad \vhat_i = \Psi(\vwake_i),
\label{eq:dim_red}
\end{equation}
where $\Psi : \mathbb R^P \to \mathbb R^L$. We specify that the latent dimension $L \ll P$ has to be selected a priori.
}

Such phase can be approached by making use of linear or non-linear techniques, such as the \podec{} and the usage of an Autoencoder (\autoe{}), respectively. 

\bmhead{Proper Orthogonal Decomposition}

In the first case, the offline part consists \RC{of} the computation of a reduced basis, composed \RC{of} a reduced number of vectors named \emph{modes}. The main assumption on which it is based is that each snapshot can be approximated by a linear combination of modes:
\[
\vwake_i  \simeq \sum_{k=1}^{\RB{L}} a_i^k \boldsymbol{\varphi}_k, \quad\quad \RB{L}\ll \RB{P}, \, i=1,\dotsc,M,
\]
with $\boldsymbol{\varphi}_i \in \mathbb{R}^P$ are the modes and $a_i^k$ are the related modal coefficients.

The computation of modes in the \podec{} procedure can be done in different ways, such as via Singular Value Decomposition (SVD) or the velocity correlation matrix. In the first case, the matrix $\mathbf{Y}$ can be decomposed via singular value decomposition in the following way:
\[\mathbf{Y} = \mathbf{U} \mathbf{\Sigma} \mathbf{V}^T, \] 
where the matrix $\mathbf{U} \in \mathbb R^{P \times \RB{L}}$ stores the \podec{} modes in its columns and matrix $\mathbf{\Sigma} \in \mathbb R^{\RB{L} \times \RB{L}}$ is a diagonal matrix including the singular values in descending order.
The matrix of modal coefficients \RB{--- so, the reduced coordinates ---} in this case can be computed by:
\[\RB{\hat{\mathbf Y}}= \mathbf{U}^T \mathbf{Y},\]
\RB{where the $\RB{\hat{\mathbf Y}} \in \mathbb R^{\RB{L}\times M}$ columns are the reduced snapshots.}
In the second case, the \RB{POD space $\mathbb{V}_{POD} = \text{span} \{[\mathbf{\phi}_i]_{i=1}^L$\} is found solving the following minimization problem:
\[\mathbb{V}_{POD} = \text{arg}\, \text{min}_{n=1,\dots, M} \dfrac{1}{M}\sum_{n=1}^{\RB{M}} \| \mathbf{Y}^n-\sum_{i=1}^{\RB{L}} a_n^i \boldsymbol{\varphi}_i \|^2 \, \forall n=1, \dots , M , \]}
where $\mathbf{Y}^n$ is the $n$-th column of the matrix $\mathbf{Y}$.
This problem is equivalent to computing the eigenvectors and the eigenvalues of the velocity correlation matrix: 
\[(\mathbf{C})_{ij}=(\mathbf{Y}^i, \mathbf{Y}^u)_{L^2(\Omega)},\] 
where $\Omega$ is the domain on which the velocity is defined (the propeller wake plane in this case).
The \podec{} modes are expressed as:
\[\boldsymbol{\varphi}_i = \dfrac{1}{\RB{M} \lambda^u_i} \sum_{j=1}^{\RB{M}} \mathbf{Y}^j V^y_{ij},\]
where $V^y$ stores the eigenvectors of the correlation matrix in its columns and $\lambda^y_i$ are its eigenvalues.

\bmhead{Autoencoder}

\RC{\autoe{} refers to a family of neural networks that, for its architectural features, has become a mathematical tool for dimensionality reduction~\cite{lee2020model}.}
In general, an \autoe{} is a neural network that is composed \RC{of} two main parts:
\RC{
\begin{itemize}
    \item the encoder: a set of layers that takes as input the high-dimensionality vector(s) and returns the reduced vector(s).
    \item the decoder, on the other hand, computes the opposite operation, returning a high-dimensional vector by passing as input the low-dimensional one.
\end{itemize}

The layers composing the \autoe{} could be in principle of any type --- e.g. convolutional~\cite{romor2022non}, dense --- but in this work}
both the encoder and the decoder are feed-forward neural networks.
\RC{For sake of simplicity, we assume here that both the encoder and the decoder are built with only one hidden layer and we denote by $\mathcal{D}$ the decoder and with $\mathcal{E}$ the encoder.} If \RB{$\vwake$} is the input of the \autoe{}, we denote by $\RB{\vtilde}=(\mathcal{D} \circ \mathcal{E})(\RB{\vwake}) = \autoe (\RB{\vwake})$ the output of the encoder, where formally:
\[
\begin{split}
&\mathcal{E}(\RB{\vwake}) = \sigma(W\RB{\vwake}+b)=\RB{\vhat},\\
&\mathcal{D}(\RB{\vhat}) = \sigma(W'\RB{\vhat}+b')=\RB{\vtilde}.
\end{split}
\]

\begin{figure}
    \centering
    \includegraphics[width=0.7\textwidth]{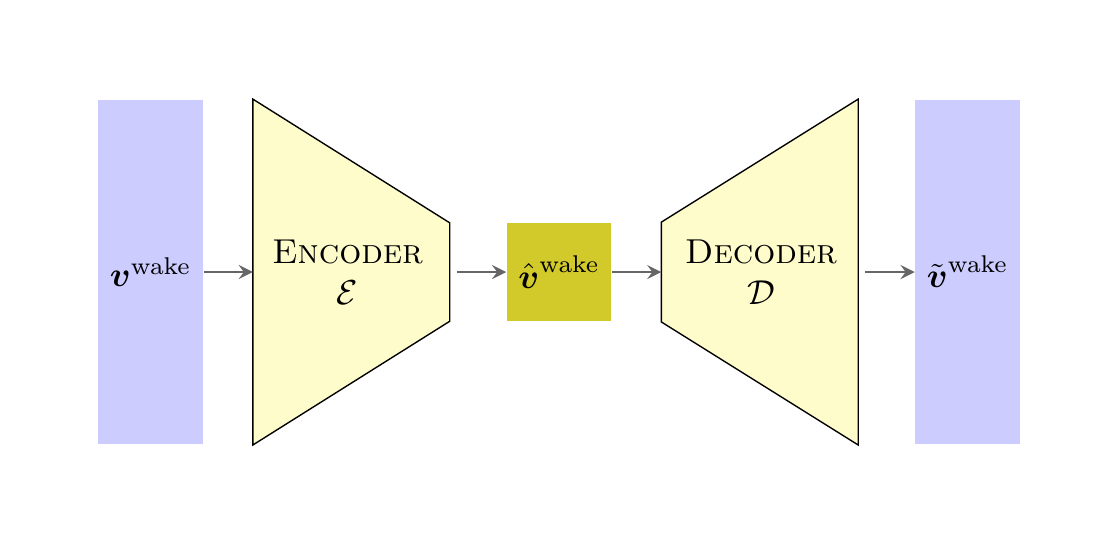}
    \caption{\RC{Schematic structure of an autoencoder.}}
    \label{fig:ae}
\end{figure}

\RC{A generic structure of an autoencoder is schematized in Figure \ref{fig:ae}.}
Weights and activation functions can be different for encoder and decoder, and of course
the case with more than one hidden layer for the encoder and the decoder can be easily derived from this formulation.
The \autoe{} is trained by minimizing the following loss function:
\[
\|\RB{\vwake}-\RB{\vtilde}\|^2=\|\RB{\vwake}-\sigma'(W'(\sigma(W\RB{\vwake}+b))+b')\|^2,
\]
\RC{where $\RB{\vwake}$ is the original (high-dimensional) vector to reduce, $\RB{\vhat}$ represents the reduced coordinates and $\RB{\vtilde}$ \RC{is} the predicted vector. In this way, the network weights are optimized such that the entire AE produces the approximation of the original vector, but compressing it onto an a priori reduced dimension. The learning step aims so to discover the latent dimension of the problem at hand.}

\RC{For what concerns the test cases here considered, two different types of autoencoders are taken into account: 

\begin{itemize}
    \item[(i)] a \emph{linear} autoencoder, i.e. without an activation function between the hidden layers, with a single hidden layer composed \RC{of} a number of neurons equal to the reduced dimension: it should exactly reproduce the behavior of the POD;
    \item[(ii)] a \emph{non-linear} autoencoder, i.e. with an activation function, and with multiple hidden layers, whose performance is compared to that of the POD.
\end{itemize}  }

\subsubsection{Approximation techniques}
The problem in the online part is to predict \RB{the (unknown) latent dimension $\vtilde$ given a new parameter $\boldsymbol{\mu}$}:
\[
\RB{\vtilde} = \pi(\boldsymbol{\mu})\quad\text{s.t.}\quad\RB{\vtilde_i} = \pi(\mupar_i) \,\text{for}\, i = 1, \dotsc, M,
\]
\RB{where $\pi : \mathbb R^p \to \mathbb R^L$ is the mapping from parameter space to reduced space. 
We can approximate such mapping by means of interpolation techniques, such as Radial Basis Functions (\rbf{}), or regression techniques, such as \ANN{}, in order to predict the latent dimension for any new parameter. Finally, the approximation of the solution in the original (high-dimensional) space requires the expansion of the reduced coordinates, which relies on the inverse compression method computed during the dimensionality reduction.}
\RA{\remark{Many approximation techniques can be employed to reconstruct the solution in reduced order models. For instance, two spread techniques are the \rbf{} and the Gaussian Process Regression (GPR). However, we remark that many other approximants can be adopted, such as the Moving Least Squares approach (MLS), which is described in detail in \cite{levin1998approximation, lancaster1981surfaces}. The reason for choosing the RBF approach is that it allows us to tune different parameters, such as the radius, and the kernel of the radial basis functions, in order to adapt our approximation to different settings of the training dataset.}}
\bmhead{Radial Basis Functions}

The \rbf{} is an approximation technique that reconstructs the original field in the following way:
\[
\RB{\vtilde} \equiv \RB{\vtilde}(\boldsymbol{\mu}) = \sum_{i=1}^M \omega_i \, \boldsymbol{\phi}(\Vert \boldsymbol{\mu}-\boldsymbol{\mu}_i \Vert),
\]
where $(\boldsymbol{\phi}(\Vert \boldsymbol{\mu}-\boldsymbol{\mu}_i \Vert))_{i=1}^M$ are called \emph{radial basis functions}, each one associated with a different center $\boldsymbol{\mu}_i$ and weighted with a weight $\omega_i$. The radial functions can have different expressions, in our case we consider the multiquadric functions $\boldsymbol{\phi}(r)=\sqrt{1+(\varepsilon r)^2}$, where $r=\Vert \boldsymbol{\mu}-\boldsymbol{\mu}_i \Vert$.

\bmhead{Artificial Neural Networks}
\RB{ 

The other technique here investigated to approximate the parametric solution manifold is \ANN{}.
The basic structure of the method is already explained in Section~\ref{subsec:ANN_param}, 

}
We consider a neural network composed \RC{of} a unique hidden layer. Its structure is:
\[\ANN{}(\boldsymbol{\mu})=\sigma(W \boldsymbol{\mu}+b)\]

The weight matrix and bias of the \ANN{} are found by training the neural network with the set of parameters and snapshots $(\boldsymbol{\mu}_i, \RB{\vwake_i})_{i=1}^M$. Then, the approximated solution $\RB{\vtilde}$ is computed from the related vector of parameters $\boldsymbol{\mu}$ as $\RB{\vtilde}=\ANN{}(\boldsymbol{\mu})$.\\

The reduced order techniques presented in this Section both for dimensionality reduction and approximation are applied in this paper to two different inverse problems. In particular, the following cases are considered:
\begin{itemize}
    \item[1.] \podec{}-\rbf{};
    \item[2.] \podec{}-\ANN{};
    \item[3.] \autoe{}-\rbf{};
    \item[4.] \autoe{}-\ANN{};
   \RALL{\item[5.] non-linear \autoe{}-\rbf{};
    \item[6.] non-linear \autoe{}-\ANN{}.}
\end{itemize}

\subsection{Wake optimization}

The construction of the reduced order model is followed by the research of the vector of parameters which better reconstructs the velocity distribution we want to reproduce. 
This investigation is addressed by solving an optimization problem, in which the aim is to minimize the difference between the approximated wake distribution predicted by the ROM and the real wake distribution.
\RA{The optimization problem can be addressed either by using a \emph{search-based}, such as the genetic algorithm (\ga{}), initially proposed in \cite{holland1973genetic}, or a \emph{gradient-based} algorithm.
In subsections \ref{opt_case1} and \ref{opt_case2} we compare the results obtained employing both approaches.
\remark{However, it is worth remarking that the genetic algorithm allows us to reach the global theoretical minimum without getting stuck into a local minimum. Gradient-based methods, instead, require derivable objective functions and get trapped in local minima in non-convex optimization.
The genetic algorithm requires a high number of evaluations, which is not a real issue since the employment of the reduced model, but it is able to converge to the global minimum. In a data-driven ROM framework, as the one proposed in this manuscript, the solution manifold is approximated with regression techniques, without any warranties on convexity. Thus, genetic methods offer a robust approach in this context, as demonstrated by its employment in similar frameworks~\cite{demoortaligustinrozzalavini2020bumi,DemoTezzeleMolaRozza2021JMSE}.}}

We dedicated this section to provid\RC{ing} a basic introduction to the genetic method, retaining a full discussion out of the topic of the present work. For a deeper focus on genetic optimization, we refer the reader to the original contribution. 
The first step of the algorithm is the definition of a population composed \RC{of} $N_{pop}$ individuals, in our case vectors of parameters $(\boldsymbol{\mu}_j)_{j=1}^{N_{pop}}$, composed \RC{of} $p$ genes, $\boldsymbol{\mu}_j \in \mathbb{R}^p$.
Then, we proceed by defining the objective function which should be minimized. We indicate by $\RB{\vtilde_j} \equiv \RB{\vtilde}(\boldsymbol{\mu}_j)$ the approximation of the wake distribution computed with the reduced order model that we are taking into account. We call $\RB{\vtarget}(\mupar)$ the real wake distribution. 

The objective function defined for each individual \RC{in} the population is:
\begin{equation}
\mathcal{F}(\boldsymbol{\mu}_j)=\Vert \RB{\vtilde}(\mupar _j)- \RB{\vtarget}(\boldsymbol{\mu}_j)\Vert .
\label{obj}
\end{equation}
The \ga{} consists in an iterative process composed \RC{of} three main steps: \emph{selection}, in which the best individuals are chosen; \emph{mate} or \emph{crossover}, where the genes of the best individuals are combined according to a certain mate probability; \emph{mutation}, changing some of the genes of the individuals.
This process is iterated a number of times which is named \emph{number of generations}.
Regarding the technical side, the \ga{} has been performed using the open-source package DEAP~\cite{DEAP}.

\section{Numerical results}
\label{sec:results}
In the present Section, the methods presented in Section \ref{sec:methods} are applied to the test case of the flow in a circular cylinder. 

\subsection{The inverse problem}
\label{fom}

\begin{figure}[ht]
\centering
\includegraphics[width=.5\textwidth]{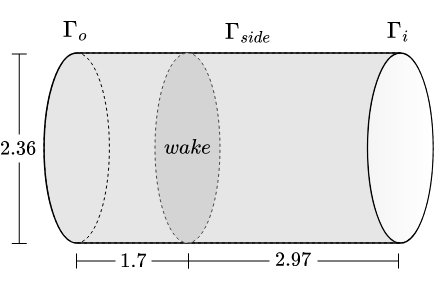}
\caption{The computational domain.}\label{fig:domain}
\end{figure}

The computational domain is a circular cylinder with height $4.67 \, m$, diameter $2.36 \, m$ and it is schematized in Figure \ref{fig:domain}, where the inlet is indicated as $\Gamma_i$, the outlet as $\Gamma_o$ \RB{and the lateral surface of the cylinder as $\Gamma_{side}$}.  The aim is to reconstruct the inlet velocity distribution given the velocity distribution at the so-called wake, which is a plane placed at a distance of $2.97$ meters from the inlet plane, as showed in \ref{fig:domain}.
This type of problem is known as \emph{inverse problem}, which is in this case applied in a \cfd{} setting. 
The physical problem at hand is modelled by the Navier-Stokes Equations (NSE) for incompressible flows. We call the fluid domain $\Omega \in \mathbb{R}^3$, $\Gamma$ its boundary; $t \in [0,T]$ is the time, $\mathbf{u}=\mathbf{u}(\mathbf{x},t)$ is the flow velocity vector field and $p=p(\mathbf{x},t)$ is the normalized pressure scalar field divided by the fluid density,
\RB{$\nu=$ \SI{1.124e-6}{\metre^2 \per \second}} is the fluid kinematic viscosity. The strong form of the NSE is the following.

\begin{subnumcases}{\label{NSE}}
\frac{\partial \mathbf{u}}{\partial t}=-\nabla \cdot (\mathbf{u} \otimes \mathbf{u})+\nabla \cdot \nu \left(\nabla \mathbf{u} + (\nabla \mathbf{u})^T \right)-\nabla p & in $\Omega \times [0,T]\, , $ \label{mom_NSE}\\
\nabla \cdot \mathbf{u} = \mathbf{0} & in $\Omega \times [0,T]\, ,$ \label{cont_NSE}\\
+ \text{ boundary conditions }& on  $\Gamma \times [0,T]\, ,$ \label{bound_NSE}\\
+ \text{ initial conditions }& in  $(\Omega,0)\, .$ \label{init_NSE}
\end{subnumcases}
\RB{The boundary and initial conditions appearing in \eqref{NSE} are the following:

\begin{itemize}
    \item at the inlet $\Gamma_i$:
    $
    \begin{cases}
    \mathbf{u}=\mathbf{u_0},\\
    \nabla p \cdot \mathbf{n} = 0,
    \end{cases}
    $\\ where $\mathbf{u_0}=\mathbf{u_0}(\mathbf{x})=(0,f(\mathbf{x}),0)$ and $f(\mathbf{x})$ is the distribution set at the inlet, specified in \eqref{dist} and \eqref{dist2}.
    \item at the outlet $\Gamma_o$:
    $
    \begin{cases}
    \nabla \mathbf{u} \cdot \mathbf{n} = \mathbf{0},\\
    p=0.
    \end{cases}
    $
    \item on the walls $\Gamma_{side}$:
    $
    \begin{cases}
    \mathbf{u} \cdot \mathbf{n} =0, \\
    \nabla p \cdot \mathbf{n} = 0.
    \end{cases}
    $
    \item moreover, $\mathbf{u}(\mathbf{x}, t) = \mathbf{0}$ and $p(\mathbf{x}, t)=0$, $\forall \mathbf{x} \in \Omega$ at initial time $t=0$.
    
\end{itemize}}

\begin{figure}
    \centering
    \includegraphics[width=0.3\textwidth, trim={30cm 0 30cm 0}, clip]{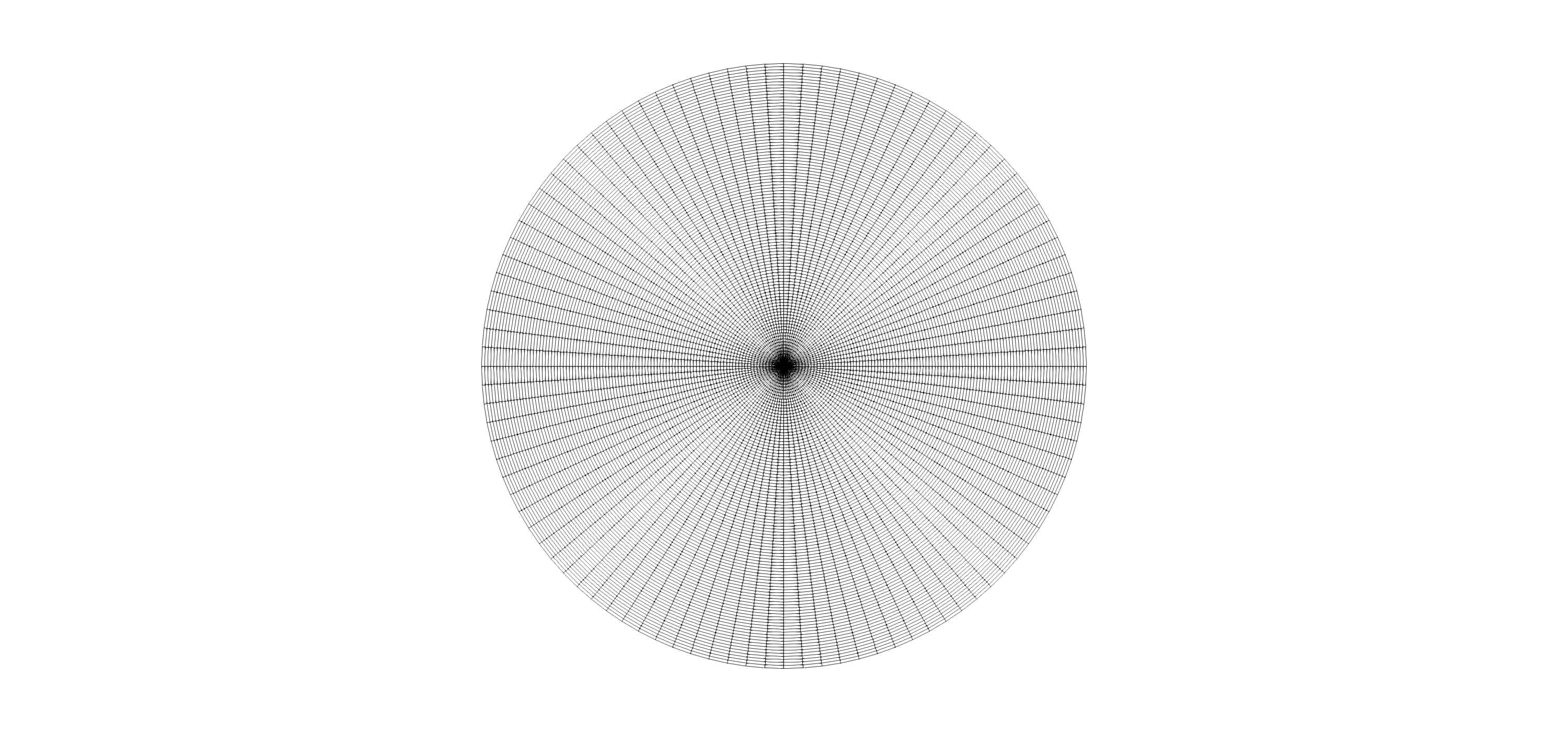}
    \includegraphics[width=0.68\textwidth]{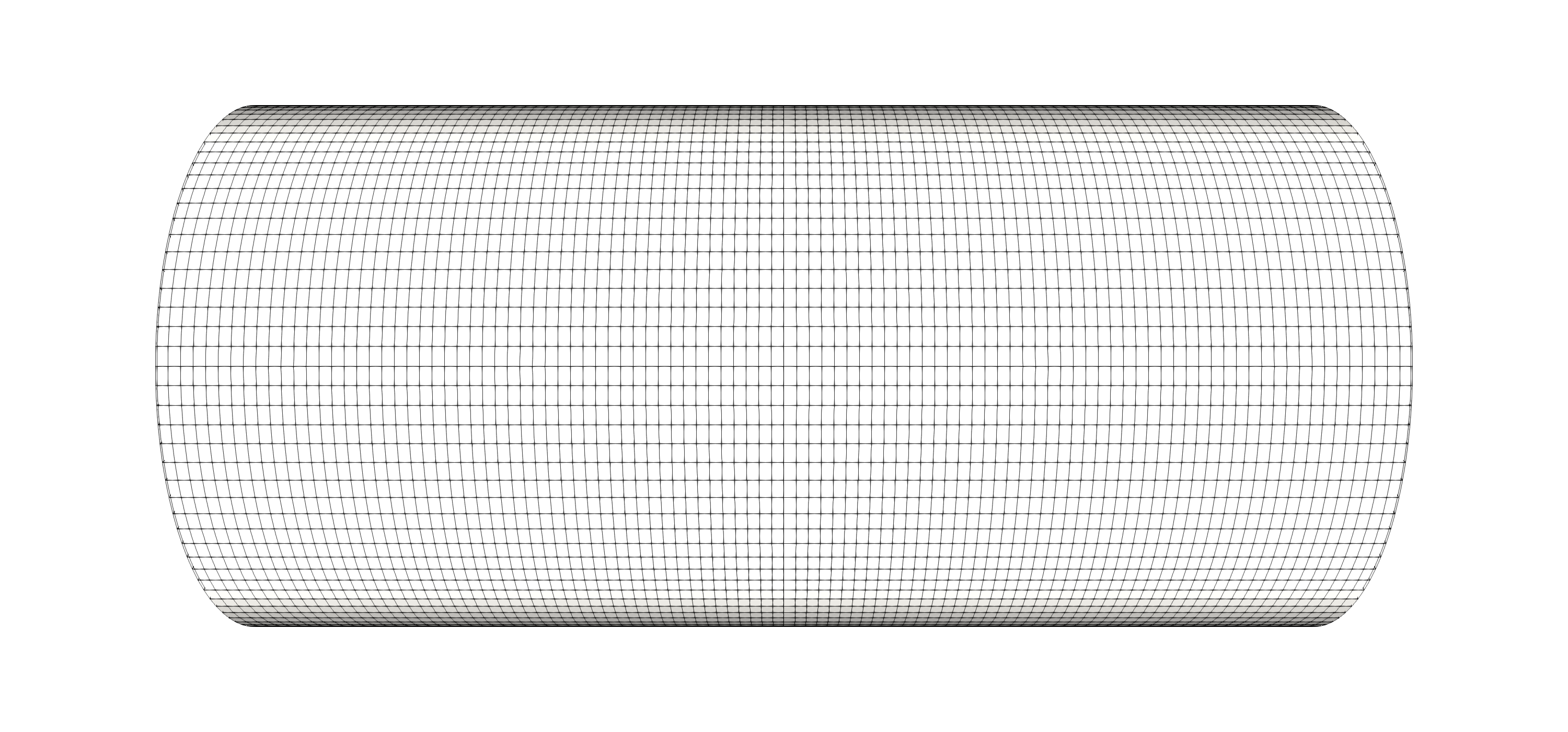}
    \caption{\RB{Representation of the mesh on a slice orthogonal to the cylinder axis (left), and on the walls (right).}}
    \label{fig: mesh}
\end{figure}

In the computation of the full order solutions of the NSE in \eqref{NSE}, the finite volume discretization is employed, by means of the open-source software OpenFOAM \cite{ofsite}. The finite volume method \cite{moukalled2016finite} is a mathematical technique that converts the partial differential equations (the NSE in our case) defined on differential volumes in algebraic equations defined on finite volumes. 
The computational mesh considered in this test case \RB{has been generated by \emph{blockMesh} and \emph{snappyHexMesh} and it is} composed by $\num{1e6}$ cells.
\RB{The mesh is a regular radial mesh, which presents 100 refinements both in the radial, the tangential and the axial direction, as can be seen from Figure \ref{fig: mesh}}.

The turbulence treatment at the full order level is characterized making use of the RANS (Reynolds Averaged Navier--Stokes) approach. This approach is based on the Reynolds decomposition, which was proposed for the first time by Osborne Reynolds \cite{reynolds1895iv}, in which each flow field is expressed as the sum of its mean and its fluctuating part.
The RANS equations are obtained by taking the time average of the NSE and adding a closure model for the well-known Reynolds stress tensor. The closure model considered in the full order model in OpenFOAM is the $\kappa-\omega$ model, proposed in its standard form in \cite{kolmogorov1941equations}, and in the SST form in \cite{menter1994two}.
This model is based on the Boussinesq hypothesis, which models the Reynolds stress tensor as proportional to the so-called \emph{eddy viscosity} and it is based on the resolution of two additional transport equation for the kinetic energy $\kappa$ and for the specific turbulent dissipation rate $\omega$.
In the full order model, we consider the SST $\kappa-\omega$ model.
\RB{The CFD simulation is run making use of the PIMPLE algorithm until convergence to a steady state.}

In this paper, the inverse problem presented is applied for two different wake distributions. The first one is a smooth function $f:\mathcal{W}\rightarrow \mathbb{R}$ defined in all the wake plane $\mathcal{W}$, written as a function of the radial coordinate:
\begin{equation}
f(r) = \sin{(4r)^2}-10.
\label{dist}
\end{equation}

The second distribution is given as a set of pointwise observations only in a selected region of the wake plane.
For the sake of simplicity only the distribution along the axis of the cylinder, which is the main direction of the flow, is taken into account.  The contributions along the secondary directions are ignored. 
\RALL{Thus, we denote with \emph{wake distribution} or \emph{velocity distribution} only the component of the velocity along the main direction of the flow.}
The observation\RC{s} are computed using the function $f$ defined as
\begin{equation}
f(x, y) = \frac{\sin{(\pi x)}\sin{((y+0.2)\pi)}}{1.2^{e^{x+y}}},
\label{dist2}
\end{equation}
where $x$ and $y$ are the cartesian coordinates in the wake plane. In total we collect $36$ observations by sampling with equispaced cartesian grid the domain $[-0.5, 0.5]^2$ in the wake plane.

\subsection{First test case: a smooth distribution}
\label{sec:smooth}
In this Section of the results, the problem considered is the reconstruction of the inlet velocity distribution when the wake velocity has the smooth distribution in \eqref{dist}.

\subsubsection{Parametrization using \nn{}}
\label{sec:NN_param_smooth}
\begin{figure}[ht]
\centering
\includegraphics[width=0.8\textwidth, ]{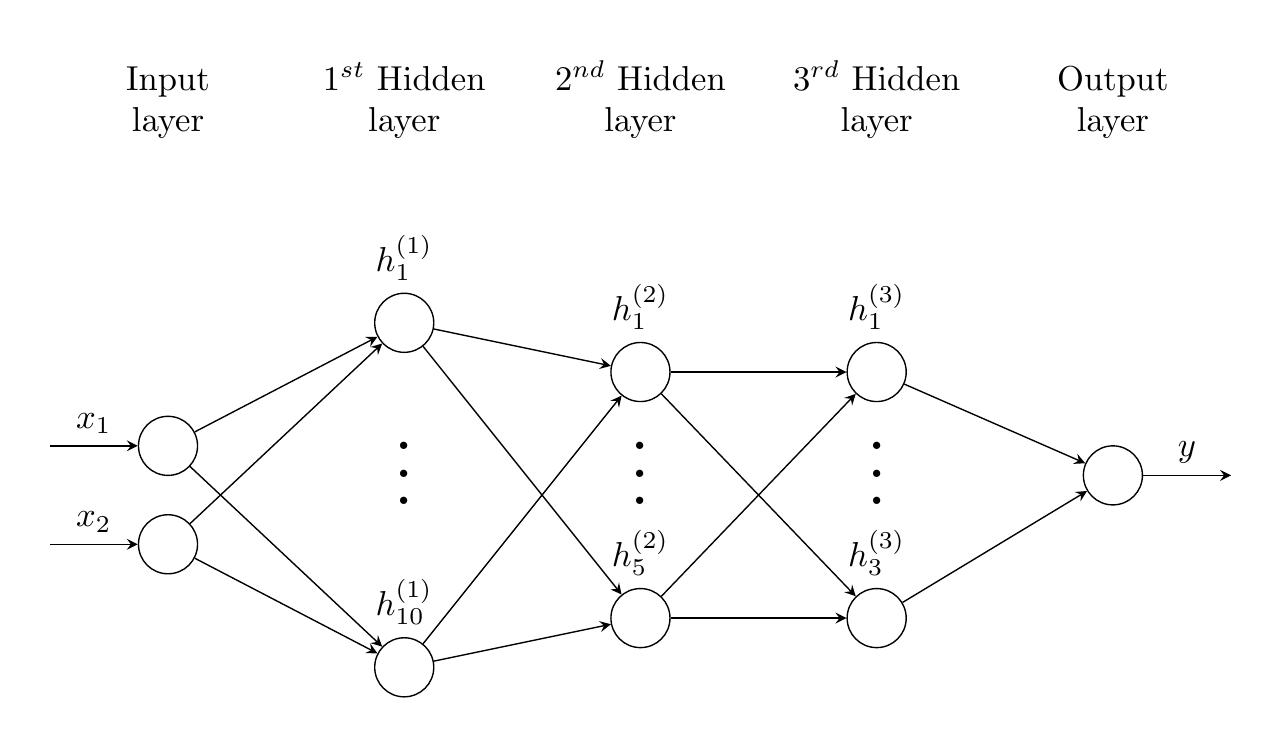}
\caption{Structure of the \nn{} used for parametrization.}\label{fig:nn1}
\end{figure}

The first step in the resolution of the inverse problem is the parametrization of the real velocity distribution at the wake through a fully-connected \nn{}. The \nn{} is represented in Figure \ref{fig:nn1} and it is composed by 3 hidden layers, with 10, 5 and 3 neurons, respectively. The inputs are the polar coordinates of the wake points and the output is the wake velocity distribution; the degree of freedom of the distribution is $P=10001$, \RB{which corresponds to the number of points at the inlet}.

As explained in Section \ref{subsec:ANN_param}, the parameters coincide with a subset of weights of the \nn{}. In particular, our choice is to consider 4 parameters, given by the biases of the last two layers.

The main idea is to generate $M$ different distributions from the \nn{} by \RB{randomly} perturbating the weights considered as parameters. Those distributions are set as inlet distributions in the full order model. The perturbated parameters and the wake distributions obtained from the full order computation are considered as parameters and snapshots, respectively, for the formulation of the \rom{}.

\begin{figure}[ht]
\centering
\includegraphics[width=\textwidth, trim={10cm 0 4cm 0}, clip]{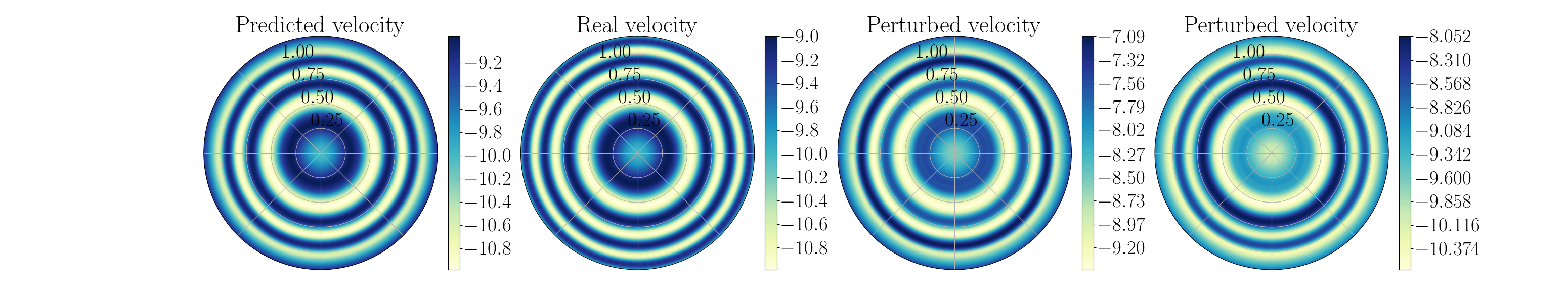}
\caption{\emph{First test case}: parametrization of the given target sinusoidal function using a \nn{}: from left to right are sketched the distribution after the training phase, the original distribution, and two perturbed configurations.}\label{fig:nn_param}
\end{figure}

In Figure \ref{fig:nn_param}, the following distributions are represented from the left to the right: the distribution predicted from the \nn{} considering the same parameters as in the training stage; the real target wake distribution; two among the $M$ distributions obtained by perturbating the parameters.

\subsubsection{Model order reduction}
\label{sec:reduction_smooth}
\RALL{The distributions obtained by the parametrization described in Section \ref{sec:NN_param_smooth} are then used as inlet initial conditions to run a set of $M$ offline simulations, following the setting described in Section \ref{fom}. The perturbed parameters of the \ANN{} and the wake distributions found from high-fidelity simulations are then used to perform a model reduction.}
In particular, in this Section different types of techniques for the reduction and approximation stages of the \rom{} are evaluated and compared and the models considered are here listed:
\begin{itemize}
    \item[1.] \podec{}-\rbf{}{};
    \item[2.] \podec{}-\nn{};
    \item[3.] \autoe{}-\rbf{};
    \item[4.] \autoe{}-\nn{}.
    \item[5.] Non linear \autoe{}-\rbf{};
    \item[6.] Non linear \autoe{}-\nn{}.
\end{itemize}

\RALL{The details related to the structure of the neural networks and the hyperparameters used for training, i.e. the learning rate, the stopping criteria, the number of hydden layers and of neurons for each layer, are defined in the discussion Section \ref{sec:discussion}.}

\begin{figure}[ht]
\centering
\subfloat[\podec{} modes]{\includegraphics[width=.9\textwidth]{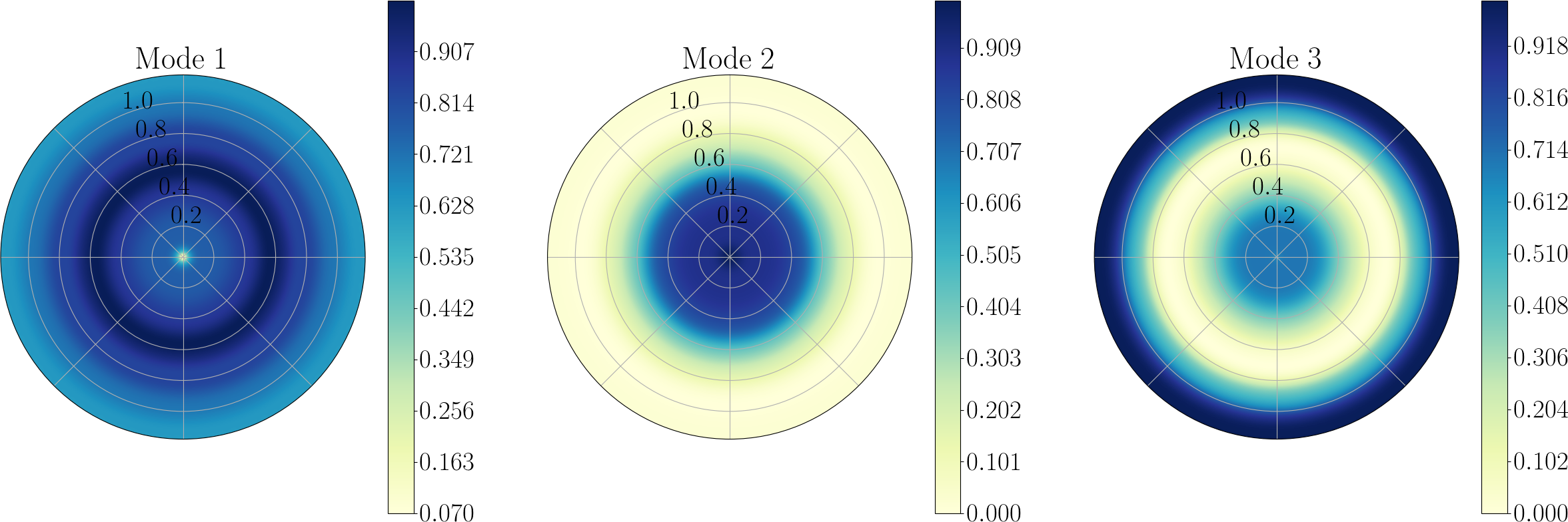}}\\
\subfloat[\autoe{} modes without activation function]{\includegraphics[width=.9\textwidth]{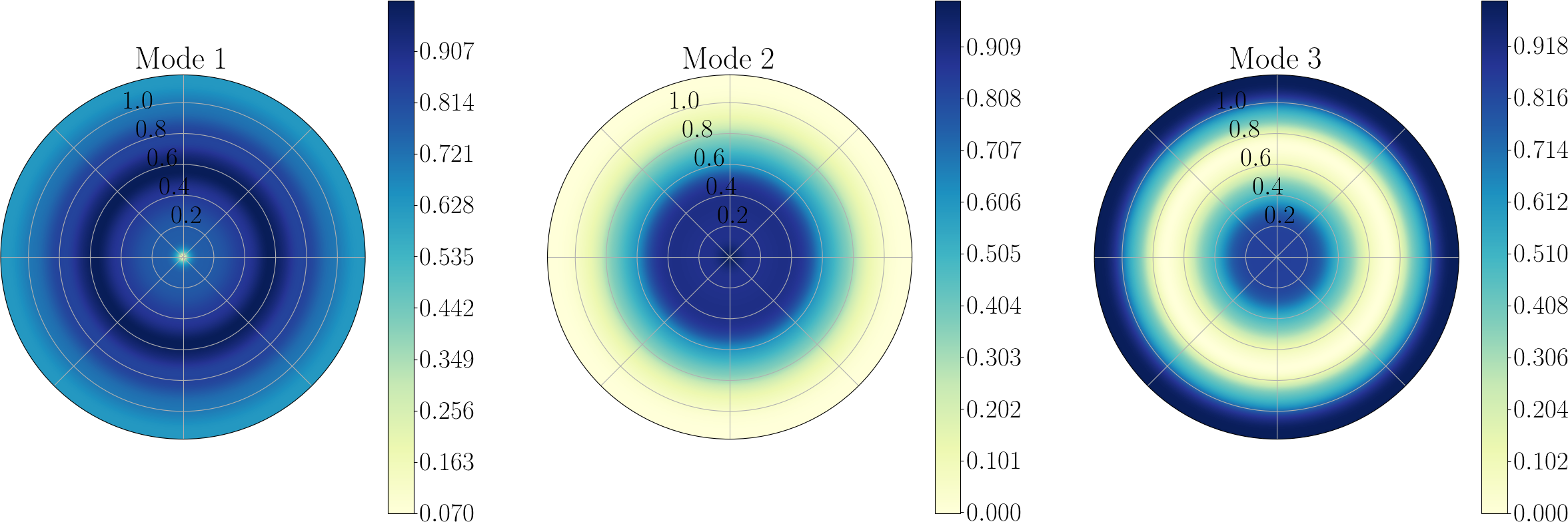}}
\caption{The \podec{} and \autoe{} modes. Modes in (b) are obtained considering a linear autoencoder.}\label{fig:modes}
\end{figure}

First of all, we consider a comparison regarding the dimensionality reduction step.
\RALL{The reduced dimension considered in this test case is $3$ and we consider a linear \autoe{}, i.e. without any activation function, composed by a single hidden layer of 3 neurons for both the encoder and the decoder parts. 
The idea is indeed to reproduce the behaviour of the \podec{}.}
A preliminary study is the graphical comparison between the \podec{} and the \autoe{} modes.
In particular, the \autoe{} modes are computed by taking as input to the decoder part of the network the identity matrix of size 3. Since the \autoe{} does not include any non linear computation, the operation is linear and ideally identical to the computation of the first 3 \podec{} modes.
From Figure \ref{fig:modes}, the shapes of the first 3 \podec{} and \autoe{} modes appear similar.

\begin{figure}[ht]
\centering
\includegraphics[width=\textwidth]{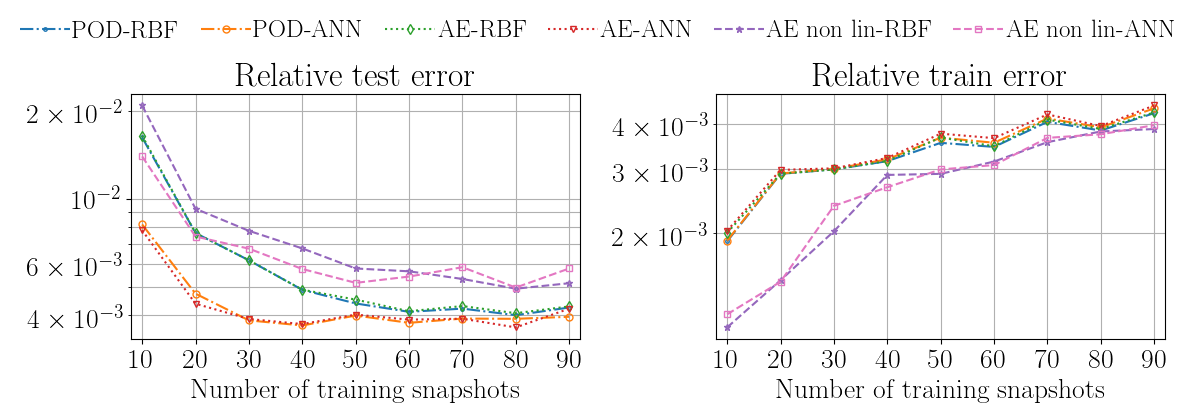}
\caption{\emph{First test case}: sensitivity analysis for all the \rom{}s with the reduced dimension $L=3$. Error is shown for test (left) and training (right) datasets.}\label{fig:reconstruction}
\end{figure}

In Figure \ref{fig:reconstruction} we provide a representation of the error \RALL{in the $L^2$ norm for: (a) a fixed test dataset; (b) the training dataset used to build each \rom{}.}
In this sensitivity analysis, the test dataset is fixed and composed by 10 snapshots and the number of snapshots used to train the reduced models varies between 10 and 90.

\RALL{Considering a dataset composed by $n$ element, we first define the relative error for each $i$-th element, say:
\[
\varepsilon_i = \dfrac{\| \vtilde_i-\vwake_i\|}{\|\vwake_i\|},
\]
where $\vtilde_i$ is the velocity distribution predicted by the reduced model starting from the test parameter, and $\vwake_i$ is the full-order snapshot in the test dataset. Then, the mean over all the elements of the test dataset is computed: $(\sum_{i=1}^{n} \varepsilon_i)/n$. 
This definition applies to both the test and the training datasets and it is used to compute the errors represented in Figure \ref{fig:reconstruction}.
The neural networks appearing in \rom{}s are trained 3 times and Figure \ref{fig:reconstruction} provides the mean test and train errors among the 3 different runs, to ensure a more reliable analysis of the \rom{}s accuracy.
}

As can be seen from Figure \ref{fig:reconstruction}, the results obtained for the test error of \podec{}-\rbf{} and \autoe{}-\rbf{} are similar to each other and, as intuitively expected, the general trend is decreasing as the number of training snapshots increases. 
The reduced methods which include a regression technique for the approximation, i.e. the \nn{}, produce better results in terms of test error with respect to the interpolating technique. \RALL{In particular, the regression technique outperforms the classical \rbf{} if we consider a small number of training snapshots; however, the results obtained by the two approximating techniques are similar when a larger training dataset is used to train the \rom{}.}


\RALL{\remark{The behavior of the neural networks (both the \autoe{} used for reduction and the \ANN{} used as approximant) strongly depends on the hyperparameters of the networks. 
This instability occurs also in the case of \emph{nonlinear \autoe{}}. Indeed, in this case the training error outperforms the results obtained with the \podec{}, but the test error remains close to the one obtained with a standard \podec{} approach, leading to the \emph{overfitting} phenomenon.
A wider discussion on the stability issues related to neural networks will be undertaken in Section \ref{sec:discussion}.}
}

\subsubsection{Wake optimization}
\label{opt_case1}

In this Section, \RA{two different optimization methods are exploited to find the combination of parameters that provides the highest accuracy in the reconstruction of the original wake. In particular, the results obtained with a genetic algorithm are compared to those obtained with a gradient-based method.} \RA{Both optimization processes exploit the \rom{} potentiality to predict the velocity distribution corresponding to each evaluation point in the parameter space.}
The reduced model taken into account is \RA{\podec{}-\nn{}, which provided the highest precision in the sensitivity analysis in \ref{sec:reduction_smooth} when the number of snapshots used for training is $90$}.

The \ga{}, composed by the selection, crossover and mutation steps, is evolved for \RA{20} generations, considering an initial population of 200 individuals, where each individual is represented by a different combination of the four parameters used for the wake parametrization in Section \ref{sec:NN_param_smooth}. \RB{All the attributes of each individual are in the interval $[-0.5, 0.5]$.}
The fitness of each individual of the population\RB{, i.e. the objective function minimized in the \ga{},} is the relative error of the wake obtained from the \rom{} with respect to the \RB{reference wake} we want to reproduce. 
\RB{We specify here more details about the implementation of the genetic evolution:
\begin{itemize}
    \item in the \emph{mate} stage, a blend crossover is performed, that modifies in-place the input individuals;
    \item in the \emph{mutation} step, a Gaussian mutation of mean $\mu=0$ and standard deviation $\sigma=1$ is applied to the input individual. The independent probability for each attribute to be mutated is set to $0.5$;
    \item the \emph{evolutionary algorithm} employed in the genetic process is the $(\mu + \lambda)$ algorithm, which selects the $\mu=50$ best individuals for the next generation and produces $\lambda=100$ children at each generation. Moreover, the probabilities that an offspring is produced by crossover and by mutation are set to $0.4$ and $0.6$, respectively.
\end{itemize}} 

\begin{figure}[ht]
\centering
\includegraphics[width=\textwidth]{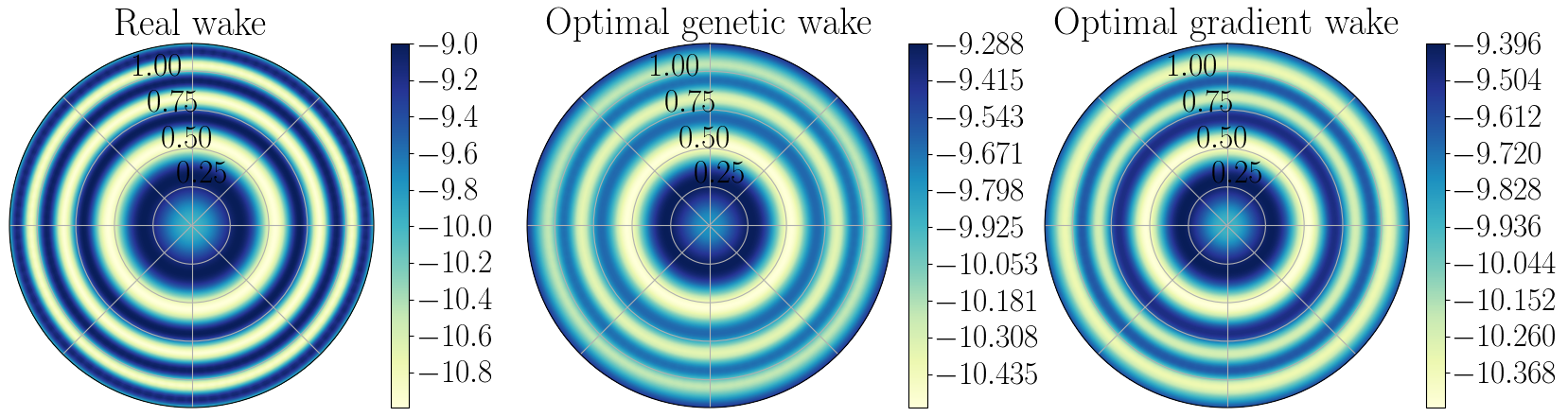}
\caption{\emph{First test case}: graphical representations of the exact benchmark wake (left), the optimal wakes obtained with a genetic (center) and gradient (right) optimization.}
\label{fig:opt}
\end{figure}

\RA{On the other hand, for the gradient optimization we considered the \emph{BFGS} (\emph{Broyden–Fletcher–Goldfarb–Shanno}) algorithm. The starting point chosen in the function evaluation is $[0.5, -0.5, -0.5, 0.5]$. The algorithm evaluated the derivative of the objective function via forward differences, setting the absolute step size to $\num{1e-17}$.

The optimal individual obtained from the genetic optimization is $\begin{bmatrix}-0.157& -0.179& -0.157& 0.0851 \end{bmatrix}$, whereas the optimal individual predicted by the gradient-based algorithm is $\begin{bmatrix} 0.00611 & -0.117 & -0.0095 &  0.0143\end{bmatrix}$.
The final value of the objective function for the optimal individual is $\sim 3.89 \, \%$ for the genetic algorithm, and $\sim 3.94 \, \%$ for the gradient algorithm.
The graphical representations of the optimal and of the real distributions are displayed in Figure \ref{fig:opt}.}

\RA{\remark{Although the two methods reach similar precision, the genetic algorithm is preferred since the gradient method is significantly influenced by the choice of the initial guess and by the step size used to evaluate the approximated derivative. We will show the numerical evidence for this argumentation in Section \ref{sec:discussion}.}}

\subsection{Second test case: a set of pointwise observations}
\label{sec:points}

In practical applications, the available data is typically provided not in the form of a smooth distribution, but in a \RC{small} number of pointwise observations, which are usually localized on a reduced part of the computational domain. 
In this section, the known data at the wake is composed of a set of 36 measurements placed inside a square of side 1 $m$ centered on the wake plane. The reference distribution considered is expressed in \eqref{dist2}.

The logical passages followed in Section \ref{sec:smooth} are reproduced in the present section with a different set of data.

\subsubsection{Parametrization using \nn{}}

\begin{figure}[ht]
\centering
\includegraphics[width=1\textwidth]{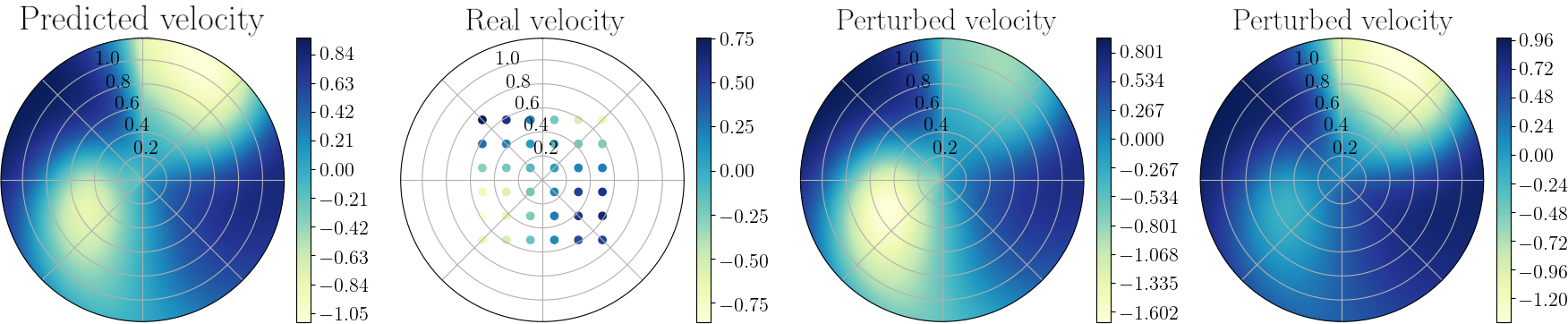}
\caption{\emph{Second test case}: parametrization of the given target  using a \nn{}: from left to right are sketched the distribution after the training phase, the original observations, and two perturbed configurations.}\label{fig:nn_param2}
\end{figure}

The \nn{} used to parametrize the wake in this test case is composed of 3 hidden layers, the first one with 8 neurons and the other two with 2 neurons. The number of parameters is 6 in this application and coincides with the perturbations of the weight matrix and the bias of the last hidden layer.
The \nn{} is trained for $50000$ epochs with a learning rate of $3\times 10^{-3}$. Since in this numerical experiment there is no radial symmetry (contrary to the previous test case) we add an additional contribution to the loss term to force the continuity between $0$ and $2\pi$ angles. Formally, the loss to minimize during the train is
\begin{equation}
L(\RB{\vtilde}, \RB{\vtarget}) = MSE(\RB{\vtilde},\RB{\vtarget}) + \lambda \| \nn{}(r_i, 0) - \nn{}(r_i, 2\pi) \|,
\end{equation}
where $MSE(\cdot, \cdot)$ is the mean square error, $\lambda$ is a weight of the additional contribution (here $\lambda = 1$) and $r_i$ for $i = 1,\dotsc,100$ are the equispaced sample points we spanned along the radius at the angles $0$ and $2\pi$. Figure~\ref{fig:nn_param2} shows indeed a continuous surface that keeps such a feature also after the perturbation of the weights.

As already done in Section \ref{sec:NN_param_smooth}, $M$ sets of parameters are generated as random perturbations in the interval $[-1,1]$, which will be used as parameters of the \rom{}.
The \nn{} perturbed with the generated parameters is used to obtain $M$ new velocity distributions, which are set as inlet distributions for the full-order model. Then, from the full-order simulations, the wake distributions corresponding to each set of parameters are obtained and considered as snapshots for the \rom{}.

\subsubsection{Model Order Reduction and Optimization}
\label{opt_case2}

From the sensitivity analysis in Section \ref{sec:reduction_smooth} the combination of the \autoe{} and the \nn{} proved to have the highest precision in the reconstruction of the wake. 
For this reason, the \autoe{}-\nn{} is the \rom{} which is chosen in this test case to train the optimization algorithm and to provide the best (approximated) wake.

\begin{figure}[ht]
\centering
\includegraphics[width=\textwidth]{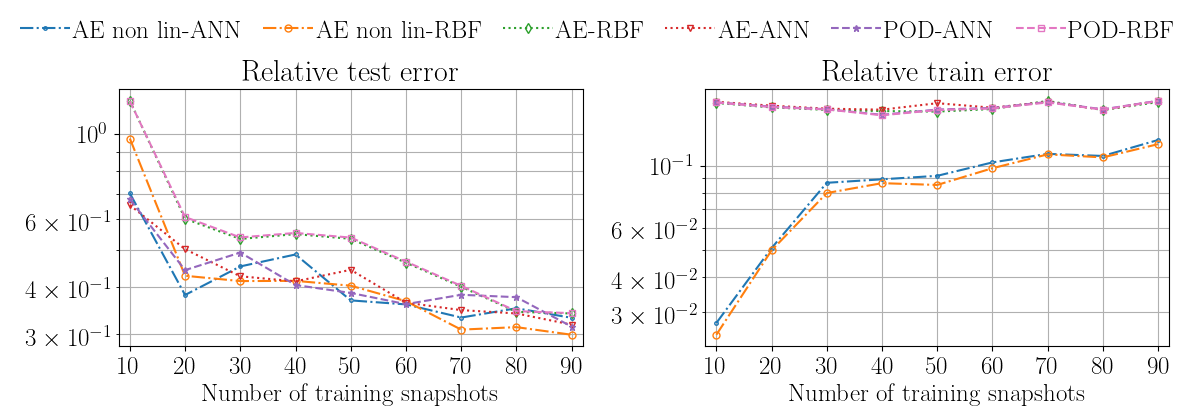}
\caption{\emph{Second test case}: sensitivity analysis for all the \rom{}s with the reduced dimension $L=4$. Error is shown for the testing (left) and training (right) datasets.}\label{fig:rom2}
\end{figure}

The \ga{}, composed by the selection, mate, and mutation, is here iterated for 20 generations. The number of individuals considered is $200$ for the first generation and $100$ for all the other generations.
The fitness is the relative error between the real wake we want to reproduce and the wake reconstructed by the \rom{}, where only the points in data are considered. 
In this experiment, the fitness decreases as the population evolves, until a percentage error of $\simeq 49 \%$ is reached. Also in this case, the \emph{BFGS} technique is not able to reach that fitness and gets blocked instead in some local minima at $61 \%$ of error with respect to the target wake. 

\begin{figure}[ht]
\centering
\includegraphics[width=.8\textwidth]{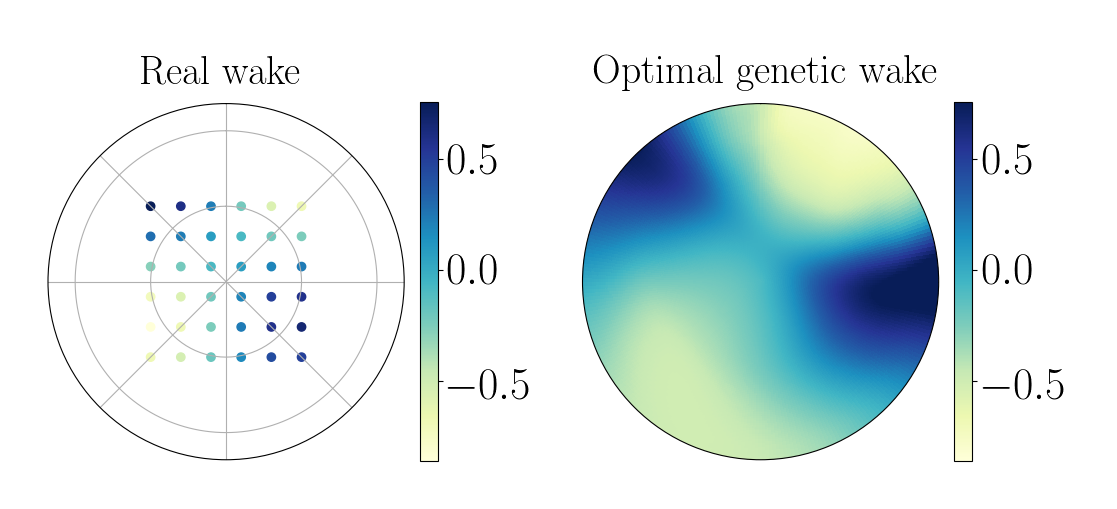}
\caption{\emph{Second test case}: the target output at wake plane (left) and the optimum representation, obtained with a genetic algorithm (right).}\label{fig:opt2}
\end{figure}

Figure~\ref{fig:opt2} shows the graphical representation of the optimum wake obtained from the \ga{}. From a graphical point of view, the optimization is able to capture the trend exhibited by the pointwise measurements. However, the error remains high after the optimization since the data in this test case is not given by a target smooth distribution, but by a few amount of points. Another possible issue is that the points are given only in a limited region of the wake, leading to a more difficult reconstruction of the distribution at the inlet, especially in the external region.

\subsection{Discussion}
\label{sec:discussion}
\RALL{
This Section is dedicated to a wider discussion of the numerical results obtained in Sections \ref{sec:smooth} and \ref{sec:points}.

Table \ref{tab:results} provides a summary of the results obtained for the two test cases. The first thing that emerges from the table is that none of the methods performs better than the others for each different value of the training snapshots. 
In general, for the first test case, the methods that provide the highest accuracy on the test snapshots are the \podec{}-\ANN{} and the linear \autoe{}-\ANN{}, as can be seen also from Figure \ref{fig:reconstruction}. On the other hand, in the second test case, the reduced order methods employing the nonlinear autoencoder provide the best results in most of the $M$ values, as can be seen in Figure \ref{fig:rom2}.
For what concerns the training error, from both Figures \ref{fig:reconstruction} and \ref{fig:rom2} we can deduce that the introduction of non-linearity in the reduction stage leads to a better reconstruction of the training snapshots, especially when $M$ is small, i.e. with a limited amount of training data.
Therefore, especially in the first test case, the autoencoder reaches a high accuracy on the train snapshots, but it is not able to generalize the information learned during training, leading to high errors on the test database. Indeed, the worst results in the first case are obtained with the nonlinear \autoe{} and are highlighted in red in Table \ref{tab:stab_first}. This behavior describes the \emph{overfitting} phenomenon, which is here mitigated by the addition of a regularization term in the loss function, say:
\[
L(\vtilde, \vtarget)=MSE(\vtilde, \vtarget) + \dfrac{\alpha}{2} \| \mathbf{W}\|^2,
\]
where $\|\mathbf{W}\|$ is the $L^2$ norm of the weight matrix, i.e. the sum over all the squared weight values, and $\alpha$ is named \emph{weight decay}.

In addition, we provide here different arguments to explain the stability issues concerning the neural networks' performance:
\begin{itemize}
    \item The networks are highly influenced by the choice of the hyper-parameters, i.e. the learning rate, the structure of the hidden layers, the stopping criteria, and the weight decay. In order to have a good result for each value of $M$, the value of the hyperparameters should be tuned \emph{ad-hoc} considering each case individually.
    In our analysis, instead, we consider the parameters fixed for each method for a fair comparison. These values are reported in Appendix \ref{sec:appendix}.
    \item The different performance of the autoencoder in the two test cases depends on the input data. We obtained better results in the second test case, where the \podec{} has poor accuracy and is not able to capture the system's behavior. To support this argumentation, we provide the results obtained in the first test case when $L=1$. From Figure \ref{fig:1mode} we can notice that the autoencoder improves the results obtained with the \podec{}
    when $M$ is low.
 \end{itemize}

\begin{table}[]
    \caption{Summary of results}
    \label{tab:results}
    \centering
    \begin{tabular}{rrcc}
    \toprule
    &  &\bf  First test case & \bf Second test case\\
    \midrule
    \multicolumn{ 2}{r}{Number of parameters $p$} & 4 & 6 \\
    \multicolumn{ 2}{r}{Reduced dimension $L$} & 3 & 4 \\
    \midrule
    \multirow{ 2}{*}{10 train snapshots} & Best method & linear AE-ANN & linear AE-ANN \\
    & Mean test error & $0.781 \, \%$ & $65.32\, \%$\\
    \midrule
    \multirow{ 2}{*}{50 train snapshots}& Best method & POD-ANN & non linear AE-ANN \\
    & Mean test error & $0.399 \, \%$ & $36.82\,\%$\\
    \midrule
    \multirow{ 2}{*}{90 train snapshots}& Best method & POD-ANN & non linear AE-RBF \\
    & Mean test error & $0.395 \, \%$ & $29.95\,\%$\\
    \midrule
    \multirow{ 2}{*}{Optimal fitness}& Genetic & $3.89 \, \%$ &  $49.83\,\%$\\
    & Gradient & $3.94 \, \%$ & $61.72\,\%$\\
    \bottomrule
    \end{tabular}
\end{table}

\begin{figure}
    \centering
    \includegraphics[width=\textwidth]{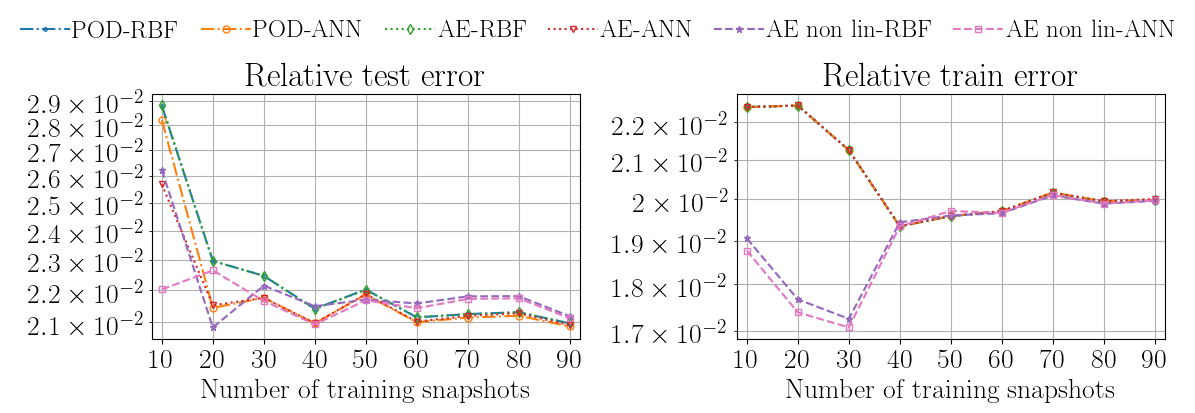}
    \caption{\emph{First test case}: sensitivity analysis for all the \rom{}s with the reduced dimension $L=1$. Error is displayed for test (left) and training (right) datasets.}
    \label{fig:1mode}
\end{figure}

\begin{table}[]
    \centering
    \caption{\emph{First test case:} stability analysis on the percentage test errors of the \rom{}s when $L=3$.}
    \label{tab:stab_first}
    \begin{tabular}{p{0.03\linewidth}%
    p{0.04\linewidth} p{0.04\linewidth}%
    p{0.04\linewidth} p{0.04\linewidth}%
    p{0.04\linewidth} p{0.04\linewidth}%
    p{0.04\linewidth} p{0.04\linewidth}%
    p{0.04\linewidth} p{0.04\linewidth}%
    p{0.04\linewidth} p{0.04\linewidth}%
    }
    \toprule
        \multirow{3}{*}{\bf M} & 
        &&
        &&
        \multicolumn{2}{c}{linear}&
        \multicolumn{2}{c}{linear}&
        \multicolumn{2}{c}{nonlinear}&
        \multicolumn{2}{c}{nonlinear}\\
        &
        \multicolumn{2}{c}{POD-RBF}&
        \multicolumn{2}{c}{POD-ANN}&
        \multicolumn{2}{c}{AE-RBF}&
        \multicolumn{2}{c}{AE-ANN}&
        \multicolumn{2}{c}{AE-RBF}&
        \multicolumn{2}{c}{AE-ANN}\\
       %
    \cmidrule(lr){2-3}
    \cmidrule(lr){4-5}
    \cmidrule(lr){6-7}
    \cmidrule(lr){8-9}
    \cmidrule(lr){10-11}
    \cmidrule(lr){12-13}
    & 
    min & max &
    min & max &
    min & max &
    min & max &
    min & max &
    min & max \\
    \midrule
    10 & 1.62 & 1.62 & 0.78 & 0.87 & 1.61 & 1.69 &  0.76 & 0.81 & 2.04 &  2.14 & 1.24 & 1.65\\
    20 & 0.76 & 0.76 & 0.42 & 0.57 & 0.76 & 0.76 &  0.43 & 0.45 & 0.9 &  0.94 & 0.64 & 0.83\\
    30 & 0.62 & 0.62 &0.38 &0.38 & 0.62 & 0.62 & 0.39 & 0.39 & 0.77 &0.78 & 0.66 & 0.69\\
    \midrule
    40 & 0.49 & 0.49 &0.37 &  0.37 & 0.49 & 0.49 &  0.37 & 0.38 & 0.63 & 0.71 & 0.57 & 0.58\\
    50 & 0.44 & 0.44 & 0.40 & 0.40 & 0.44 & 0.48 & 0.40 & 0.40 & 0.54 & 0.61 & 0.46 & 0.60\\
    60 & 0.41 & 0.41 & 0.38 & 0.38 & 0.41 & 0.42 &  0.38 & 0.41 & 0.55 & 0.60 & 0.49 & 0.59\\ 
    \midrule
    70 & 0.42 & 0.42 &0.36 & 0.44 & 0.42 & 0.43 & 0.37 & 0.41 & 0.51 & 0.57 & 0.54 &  0.63\\
    80 & 0.40 & 0.40 & 0.36 & 0.41 & 0.40 & 0.42 & 0.36 & 0.37 & 0.46 & 0.53 & 0.45 &  0.53\\
    90 & 0.43 & 0.43 &0.39 & 0.41 & 0.43 & 0.43 & 0.41 & 0.45 & 0.49 & 0.54 & 0.53 & 0.61\\
    \bottomrule
    \end{tabular}
\end{table}

\begin{table}[]
    \centering
    \caption{\emph{Second test case}: stability analysis on the percentage test errors of the \rom{}s, when $L=4$.}
    \label{tab:stab_second}
    \begin{tabular}{p{0.03\linewidth}%
    p{0.04\linewidth} p{0.04\linewidth}%
    p{0.04\linewidth} p{0.04\linewidth}%
    p{0.04\linewidth} p{0.04\linewidth}%
    p{0.04\linewidth} p{0.04\linewidth}%
    p{0.04\linewidth} p{0.04\linewidth}%
    p{0.04\linewidth} p{0.04\linewidth}%
    }
    \toprule
        \multirow{3}{*}{\bf M} & 
        &&
        &&
        \multicolumn{2}{c}{linear}&
        \multicolumn{2}{c}{linear}&
        \multicolumn{2}{c}{nonlinear}&
        \multicolumn{2}{c}{nonlinear}\\
        &
        \multicolumn{2}{c}{POD-RBF}&
        \multicolumn{2}{c}{POD-ANN}&
        \multicolumn{2}{c}{AE-RBF}&
        \multicolumn{2}{c}{AE-ANN}&
        \multicolumn{2}{c}{AE-RBF}&
        \multicolumn{2}{c}{AE-ANN}\\
       %
    \cmidrule(lr){2-3}
    \cmidrule(lr){4-5}
    \cmidrule(lr){6-7}
    \cmidrule(lr){8-9}
    \cmidrule(lr){10-11}
    \cmidrule(lr){12-13}
    & 
    min & max &
    min & max &
    min & max &
    min & max &
    min & max &
    min & max \\
    \midrule
10  & 122.08 & 122.08 & 56.46 & 80.45 & 122.34 & 122.38 & 58.16 & 75.26 & 92.35 & 101.82 & 54.26 & 100.21 \\
20  & 60.80 & 60.80 & 40.72 & 47.85 & 60.13 & 60.24 & 42.16 & 55.72 & 41.07 & 44.48 & 30.32 & 46.01 \\
30  & 53.78 & 53.78 & 41.54 & 52.91 & 53.29 & 53.30 & 38.26 & 48.82 & 40.52 & 41.93 & 37.53 & 49.45 \\
    \midrule
40  & 55.25 & 55.25 & 37.16 & 43.01 & 54.76 & 55.04 & 38.42 & 46.34 & 38.96 & 44.21 & 45.40 & 52.15 \\
50  & 53.73 & 53.73 & 35.96 & 41.43 & 53.30 & 53.37 & 40.07 & 47.69 & 36.77 & 42.99 & 32.15 & 40.05 \\
60  & 46.49 & 46.49 & 33.44 & 38.89 & 46.17 & 46.22 & 33.94 & 38.85 & 35.91 & 37.16 & 31.84 & 42.95 \\
    \midrule
70  & 40.20 & 40.20 & 37.17 & 39.85 & 39.93 & 40.00 & 30.06 & 37.57 & 28.92 & 33.38 & 25.91 & 39.48 \\
80  & 34.52 & 34.52 & 33.99 & 43.17 & 34.35 & 34.39 & 33.56 & 34.85 & 30.09 & 32.61 & 29.09 & 43.84 \\
90  & 34.10 & 34.10 & 28.99 & 35.16 & 34.03 & 34.04 & 29.24 & 35.51 & 28.60 & 30.67 & 27.26 & 37.48 \\

    \bottomrule
    \end{tabular}
\end{table}
}
\RA{An additional point that can be discussed is the type of optimization algorithm used. As pointed out in Section \ref{opt_case1}, the genetic algorithm is preferred to a gradient-based algorithm, because the gradient method can get stuck into local minima. The evidence of this statement is provided in Table \ref{tab:grad}, where we show the results obtained with the \emph{BFGS}
 method starting from different initial guesses. In the first three cases analyzed, the method converges to three different local minima; in the last two cases, it gets stuck in the first evaluation point.

\begin{table}[]
    \centering
    \caption{\emph{First test case}: results of gradient optimization for different initial guesses.}
    \label{tab:grad}
    \begin{tabular}{p{0.03\linewidth}
    p{0.13\linewidth}
    p{0.13\linewidth}
    p{0.08\linewidth}
    p{0.45\linewidth}
    }
    \toprule
         N. test & Function evaluations & Gradient evaluations & Final fitness ($\%$) & Final point   \\
         \midrule
         1 & 85 & 17 & 3.94 & $\begin{bmatrix} 0.01 & -0.12 & -0.01 & 0.01\end{bmatrix}$\\
         2 & 236 & 47 & 3.98 & $\begin{bmatrix}-0.1 & -0.06 & -0.18 &  0.1\end{bmatrix}$\\
         3 & 252 & 48 & 4.03 & $\begin{bmatrix}-0.33 & 0.03 & -0.1 &  0.1\end{bmatrix}$\\
         4 & 5 & 1 & 4.25 & $\begin{bmatrix}-0.1 & 0.1 & -0.1 &  0.1 \end{bmatrix}$\\
         5 & 5 & 1 & 5.36 & $\begin{bmatrix}0.1 & 0.1 & 0.1 &  0.1\end{bmatrix}$\\
         \bottomrule
    \end{tabular}
\end{table}

\begin{table}[]
    \centering
    \caption{\emph{Second test case}: results of gradient optimization for different initial guesses.}
    \label{tab:grad2}
    \begin{tabular}{p{0.03\linewidth}
    p{0.13\linewidth}
    p{0.13\linewidth}
    p{0.08\linewidth}
    p{0.45\linewidth}
    }
    \toprule
         N. test & Function evaluations & Gradient evaluations & Final fitness ($\%$) & Final point   \\
         \midrule
         1 & 467 & 65 & 110.30 & $\begin{bmatrix} -0.047 & -0.045 & -0.23 & -0.27 &   0.108 & 0.55\end{bmatrix}$\\
         2 & 536 & 75 & 61.72 & $\begin{bmatrix}0.27 & -0.54 & 0.107 &  0.44 & 0.35 &-0.53\end{bmatrix}$\\
         3 & 512 & 72 & 69.42 & $\begin{bmatrix}0.25 & 0.062 & -0.39 &  0.17 & -0.42 &0.33\end{bmatrix}$\\
         4 & 719 & 101 & 74.8 & $\begin{bmatrix}0.03 & -0.11 & -0.097 &  0.25 & 0.49 &-0.27\end{bmatrix}$\\
         5 & 382 & 53 & 79.38 & $\begin{bmatrix}-0.096 & 0.46 & 0.28 &  -0.027 & 0.058 &-0.49\end{bmatrix}$\\
         \bottomrule
    \end{tabular}
\end{table}
}

\section{Conclusions}
\label{sec:conclusions}
We presented in this paper a data-driven computational pipeline that allows us to efficiently face inverse problems by means of the model order reduction technique.
Neural networks are involved in the such framework at three different levels:
\begin{itemize}
    \item the parametrization of a given (generic) wake distribution by perturbating a subset of the weights of a fully-connected neural network. The parameters are exploited to obtain the inlet distributions of the full-order simulations. The wake distributions obtained from the full-order simulations are considered snapshots associated with the input parameters, obtaining a parametrization of the original problem.
    \item the compression of the dimensionality of the full-order snapshots, which typically belong to high-dimensional spaces, obtained with a properly structured autoencoder;
    \item the approximation of the solution manifold, allowing for predicting new solutions for any unseen parameters. This last task is obtained with a second fully-connected neural network.
\end{itemize}
The results obtained by both the test cases in this work are promising. The neural network parametrization shows the possibility to produce several different distributions playing with very few parameters (4 and 6), becoming a valuable tool also for such an objective. Of course, the benefits of such parametrization are strictly related to the problem at hand, since here we are assuming that the target distribution is somehow related to the boundary condition. The investigation of the parametric configuration shows also a large range of admissible distributions. However, we highlight there is no possibility to precisely select the region of interest for the parametrization and the weights of the neural network which have been perturbated have no physical significance linked to the wake velocity distribution.
Dimensionality compression and manifold approximation are already investigated uses of neural networks, and also in this contribution, the results confirm their better performances with respect to the more consolidated techniques.

Possible future developments of the presented pipeline may interest the parametrization through neural networks: at the current stage, it results in a sort of \emph{black-box} procedure, since the parameters and their ranges are selected following a trial and error process. Such an issue can be alleviated by looking at the intermediate output --- the output of the hidden layers of the network --- to propose a meaningful criterion for parameter selection.

\subsection*{Appendix}
\label{sec:appendix}
In this section, we report the specifications of the hyperparameters of the neural networks used to build the \rom{}s in the first and second test cases (Tables \ref{tab:networks_first} and \ref{tab:networks_second}).

\begin{table}[h!]
    \centering
 \caption{\emph{First test case:} neural networks setting in \rom{}s.}
    \label{tab:networks_first}
    \begin{tabular}{>{\centering\arraybackslash}p{0.18\linewidth}
    >{\centering\arraybackslash}p{0.08\linewidth}
    >{\centering\arraybackslash}p{0.12\linewidth}
    >{\centering\arraybackslash}p{0.14\linewidth}
    >{\centering\arraybackslash}p{0.06\linewidth}
    >{\centering\arraybackslash}p{0.06\linewidth}
    >{\centering\arraybackslash}p{0.12\linewidth}
    }
    \bottomrule
    \multirow{3}{*}{\textbf{NN}}& \multirow{3}{*}{Structure} & \multirow{3}{*}{Non-linearity} &\multirow{3}{*}{Learning rate} & \multicolumn{2}{c}{\multirow{2}{*}{Stop criteria}}& \multirow{3}{*}{Weight decay}\\
    &&&&&\\
    \cline{5-6}
    &&&& epochs & final loss &\\
    \midrule
\textbf{\ANN{}} &\multirow{2}{*}{$[4,4]$}&\multirow{2}{*}{Softplus}& \multirow{2}{*}{$\num{5e-3}$} & \multirow{2}{*}{200000} & \multirow{2}{*}{$\num{1e-3}$} & \multirow{2}{*}{0}\\ 
(\podec{}-\ANN{}) &&&&&&\\
\midrule
   \textbf{\ANN{}}&\multirow{2}{*}{$[4,4]$}&\multirow{2}{*}{Softplus}& \multirow{2}{*}{$\num{5e-3}$} & \multirow{2}{*}{100000} & \multirow{2}{*}{$\num{1e-4}$} & \multirow{2}{*}{0}\\
    (lin \autoe{}-\ANN{})&&&&&& \\
    \midrule
   \textbf{\ANN{}}  &\multirow{2}{*}{$[40,40]$}&\multirow{2}{*}{Softplus}& \multirow{2}{*}{$\num{5e-3}$} & \multirow{2}{*}{100000} & \multirow{2}{*}{$\num{1e-5}$} & \multirow{2}{*}{0}\\
   (non lin \autoe{}-\ANN{})&&&&&&\\
   \midrule
  \textbf{\autoe{}} &\multirow{2}{*}{$[3,3]$}&\multirow{2}{*}{None}& \multirow{2}{*}{$\num{1e-3}$} & \multirow{2}{*}{1000} & \multirow{2}{*}{} & \multirow{2}{*}{0}\\
   (linear) &&&&&&\\ 
   \midrule
  \textbf{\autoe{}} &$[200, 3,$&\multirow{2}{*}{Leaky ReLU}& \multirow{2}{*}{$\num{5e-4}$} & \multirow{2}{*}{5000} & \multirow{2}{*}{$\num{5e-4}$} & \multirow{2}{*}{$\num{5e-4}$}\\
  (non linear)&$3, 200]$&&&&&\\
  
    \bottomrule
    \end{tabular}
\end{table}

\begin{table}[h!]
    \centering
 \caption{\emph{Second test case:} neural networks setting in \rom{}s.}
    \label{tab:networks_second}
    \begin{tabular}{>{\centering\arraybackslash}p{0.18\linewidth}
    >{\centering\arraybackslash}p{0.08\linewidth}
    >{\centering\arraybackslash}p{0.12\linewidth}
    >{\centering\arraybackslash}p{0.14\linewidth}
    >{\centering\arraybackslash}p{0.06\linewidth}
    >{\centering\arraybackslash}p{0.06\linewidth}
    >{\centering\arraybackslash}p{0.12\linewidth}
    }
    \bottomrule
    \multirow{3}{*}{\textbf{NN}}& \multirow{3}{*}{Structure} & \multirow{3}{*}{Non-linearity} &\multirow{3}{*}{Learning rate} & \multicolumn{2}{c}{\multirow{2}{*}{Stop criteria}}& \multirow{3}{*}{Weight decay}\\
    &&&&&\\
    \cline{5-6}
    &&&& epochs & final loss &\\
    \midrule
\textbf{\ANN{}} & $[40,20,10]$& Softplus & $\num{2e-3}$ & &$0.1$ & 0\\
\midrule
  \textbf{\autoe{}} &\multirow{2}{*}{$[4,4]$}&\multirow{2}{*}{None}& \multirow{2}{*}{$\num{3e-3}$} & \multirow{2}{*}{5000} & \multirow{2}{*}{} & \multirow{2}{*}{0}\\
   (linear) &&&&&&\\ 
   \midrule
  \textbf{\autoe{}} &$[200, 4,$&\multirow{2}{*}{Leaky ReLU}& \multirow{2}{*}{$\num{3e-4}$} & \multirow{2}{*}{2000} & \multirow{2}{*}{} & \multirow{2}{*}{$\num{5e-4}$}\\
  (non linear)&$4, 200]$&&&&&\\
  
    \bottomrule
    \end{tabular}
\end{table}

\newpage

\section*{Acknowledgements}
This work was partially funded by INdAM-GNCS 2020-2021 projects, by European High-Performance Computing Joint Undertaking project Eflows4HPC GA N. 955558, by PRIN "Numerical Analysis for Full and Reduced
Order Methods for Partial Differential Equations" (NA-FROM-PDEs) project
by European Union Funding for
Research and Innovation --- Horizon 2020 Program --- in the framework
of European Research Council Executive Agency: H2020 ERC CoG 2015
AROMA-CFD project 681447 ``Advanced Reduced Order Methods with
Applications in Computational Fluid Dynamics'' P.I. Professor Gianluigi Rozza.

\bibliographystyle{abbrv}
\bibliography{biblio}

\end{document}